\documentclass[a4paper,11pt]{article}

\usepackage[margin=2.5cm]{geometry} 

\usepackage[utf8]{inputenc}
\usepackage[T1]{fontenc}

\usepackage{amsmath,amsfonts,amsthm,amssymb}
\usepackage{nicefrac,xfrac}
\usepackage{enumerate}
\usepackage{graphicx}
\usepackage{xcolor}
\usepackage{comment}
\usepackage{hyperref}

\DeclareMathOperator{\dist}{dist}
\DeclareMathOperator{\diam}{diam}

\newtheorem{theorem}{Theorem}[section]
\newtheorem{lemma}[theorem]{Lemma}
\newtheorem{proposition}[theorem]{Proposition}
\newtheorem{corollary}[theorem]{Corollary}

\newtheorem{problem}[theorem]{Problem}
\newtheorem{remark}[theorem]{Remark}

\title{\textbf{Exact distance Kneser graphs}}

\author{
  Agustina Victoria Ledezma\footnote{Instituto de Matem\'atica Aplicada San Luis, Universidad Nacional de San Luis and CONICET, San Luis, Argentina. e-mail: avledezma@unsl.edu.ar} 
  \and 
  Adri\'an Pastine\footnote{Departamento de Matem\'atica, Universidad Nacional de San Luis, San Luis, Argentina. e-mail: agpastine@gmail.com}
  \and 
  Mario Valencia-Pabon\footnote{Universit\'e de Lorraine, LORIA, Nancy, France. e-mail: mario.valencia@loria.fr}
}

\date{}

\begin{document}
\maketitle

\begin{abstract}
For any graph $G = (V,E)$ and positive integer $d$, the \textit{exact distance-$d$} graph $G_{=d}$ is the graph with vertex set $V$, where two vertices are adjacent if and only if the distance between them in $G$ is $d$. We study the exact distance-$d$ Kneser graphs. For these graphs, we characterize the adjacency of vertices in terms of the cardinality of the intersection between them. We present formulas describing the distance between any pair of vertices and we compute the diameter of these graphs.

\end{abstract}


\section{Introduction and main results}
Let $G$ be a connected graph. Given  two vertices $a$ and $b$ in $G$, $\dist_G(a,b)$,
 the {\em distance} between $a$ and $b$, is defined as the length of the shortest path in $G$ joining $a$ to $b$. The diameter of $G$, denoted $\diam(G)$, is the maximum distance between any pair of vertices in $G$.
 
The concept of exact distance-$d$ graph, where $d$ is a positive integer, was introduced by Ne\u set\u ril and Ossona de Mendez \cite{NeOs12} (Section 11.9). Formally, if $G=(V,E)$ is a graph, the {exact distance}-$d$ graph $G_{=d}$ of $G$ is the graph with vertex set $V$ and where two vertices $a$ and $b$ in $G_{=d}$ are adjacent if and only if $\dist_G(a,b)=d$. Note that $G_{=1} = G$. This concept is related to that of power graphs. Given a positive integer $p$, the $p$-th \textit{power graph} $G^p=(V(G), E^p)$ of $G$ is defined as the graph where the vertex set is $V(G)$ as its vertex set and the edge set $E^p$ contains an edge ${xy}$ if and only if $\dist_G(x, y) \leq p$.

The main focus in early research of exact distance graphs was on their chromatic number \cite{BEHJ19, HKQ18, NeOs12, NeOs15}. The exact distance graphs have been much earlier considered for hypercube graphs in the frame of the so-called cube-like graphs \cite{DHLL90, LMT88, Zi01}. More recently, Bre\u sar et al. \cite{BGKT19} considered the structure of exact distance graphs of some graph products. They remarked that the exact distance graphs are not only interesting because of the chromatic number, but also as a general metric graph theory concept. 
 In particular, closed formulas for the distances between pairs of vertices and for  the diameter of a graph are of interest in the study of chemical indices. In particular, the Wiener index, $W(G)$ defined as 
\[
W(G)=\sum_{\{x,y\}\in E(G)}\dist_G(x,y),
\]
has been greatly studied (see \cite{balakrishnan2007wiener,dobrynin2001wiener,dobrynin2002wiener,gutman1986mathematical,eu2006generalised,wiener1947structural,xu2006catacondensed}). In \cite{balak} formulas for the Wiener index of the Kneser graph and of the Odd graph  have been obtained.

Let $k$ and $r$ be positive integers and let $[2k+r] =
\{1,2,\ldots,2k+r\}$. Let $[2k+r]^k$ be the set of $k$-subsets of
$[2k+r]$. The Kneser graph $K(2k+r,k)$ is the graph with vertex set
$[2k+r]^k$ and where two $k$-subsets $A,B \in [2k+r]^k$ are joined by
an edge if $A \cap B = \emptyset$. Note that $K(5,2)$ is the well-known
Petersen graph. It is easy to show that the Kneser graph $K(2k+r,k)$
is a connected regular graph having $\binom{2k+r}{k}$ vertices of
degree $\binom{k+r}{k}$. Theoretical properties of Kneser graphs have
been studied in past years (see for example \cite{FUR, LOV, STA76, VV05} and references therein).

Stahl shows in \cite{STA76} the following result.

\begin{lemma}[\cite{STA76}]
Let $A,B \in [2k+r]^k$ be two vertices in $K(2k+r,k)$ joined by a
path of length $2p$ ($p \geq 0$). Then $|A \cap B| \geq k-rp$. 
\label{lstahl}
\end{lemma}
The following result is a consequence of Lemma \ref{lstahl}.

\begin{corollary}
Let $A,B \in [2k+r]^k$ be two vertices in $K(2k+r,k)$ joined by a
path of length $2p+1$ ($p \geq 0$). Then $|A \cap B| \leq rp$. 
\label{corostahl}
\end{corollary}

Valencia-Pabon and Vera \cite{VV05} studied the diameter of Kneser graphs, and obtained the following results.

\begin{proposition}[\cite{VV05}]
Let $k$ and $r$ be positive integers, with $k \geq 2$ and $r \geq k-1$. Then $\diam(K(2k+r,k))=2$.
\label{propvv}
\end{proposition}

\begin{lemma}[\cite{VV05}]
\label{lvv}
Let $A,B \in [2k+r]^k$ be two different vertices of the Kneser graph
$K(2k+r,k)$, where $1 \leq r < k-1$. If $|A \cap B| = s$ then
\[\dist(A,B) = \min\left\{2\left\lceil \frac{k-s}{r} \right\rceil, 2\left\lceil\frac{s}{r} \right\rceil + 1\right\}.\]
\end{lemma}

\begin{theorem}[\cite{VV05}]
\label{thvv}
Let $k$ and $r$ be positive integers. Then 
$$\diam(K(2k+r,k))= \left\lceil \frac{k-1}{r} \right\rceil + 1.$$
\end{theorem}

There are some generalizations to Kneser-like graphs in the literature (see for example \cite{AAC18,CheW08,Fr85,FUR} and references therein). Let $K(n,k,i)$ be the {\em generalized Kneser graph}, i.e. the graph with vertex set $[n]^k$ and the edges connecting all pairs of vertices with an intersection smaller or equal to $i$. Denote by $J(n,k,i)$ the {\em generalized Johnson graph}, i.e. the graph with the same vertex set $[n]^k$ and edges connecting all pairs of vertices with intersections exactly $i$. In particular, in \cite{CheW08} the distance between two vertices of $K(n,k,i)$ was characterized in terms of the size of their intersection, and the diameter of the graph was computed. 
These results are given in Theorems \ref{thdistgenkn} and \ref{thdiamgenkn}.

\begin{theorem}[\cite{CheW08}]
\label{thdistgenkn}
Let $A$ and $B$ be two vertices in $K(n,k,i)$, where $n=2k-i+r$ and $0\leq r<k-2i-1$. If  $|A \cap B|=s>i$ then 
\[\dist_{K(n,k,i)}(A,B)
=\min \left\{ 2 \left\lceil \frac{k-s}{i+r} \right\rceil , 2 \left\lceil \frac{s-i}{i+r} \right\rceil + 1  \right\} .
\]
\end{theorem}
\begin{theorem}[\cite{CheW08}]
\label{thdiamgenkn}
If $k$ and $r$ are positive integers, $i<k$ is a non-negative integer and $n=2k-i+r$ then 
$$\diam(K(n,k,i))= \left\lceil \frac{k-i-1}{i+r} \right\rceil + 1.$$
\end{theorem}
Later, in \cite{AAC18}, similar results were obtained for $J(n,k,i)$, which are given in Theorems \ref{thdistgenjohn} and \ref{thdiamgenjohn}.
\begin{theorem}[\cite{AAC18}]
\label{thdistgenjohn}
Let $A$ and $B$ be two vertices in $J(n,k,i)$, where $n>k>i$ are non-negative integers, 
and $(n,k,i)\neq (2k,k,0)$.
 If $|A \cap B|=s$
then 
\[
\dist_{J(n,k,i)}(A,B) = \begin{cases}
    \makebox[9cm][s]{$\displaystyle 3$ \hfill if $s < \min \{i, -n+3k-2i \}$;} \\[1.5ex]
    \makebox[9cm][s]{$\displaystyle \left\lceil \frac{k-s}{k-i} \right\rceil$ \hfill if $-n+3k-2i \leq s < i$;} \\[1.5ex]
    \makebox[9cm][s]{$\displaystyle \min \left\{ 2 \left\lceil \frac{k-s}{n-2k+2i} \right\rceil , 2 \left\lceil \frac{s-i}{n-2k+2i} \right\rceil + 1 \right\}$ \hfill if $s \geq i$.}
\end{cases}
\]
\end{theorem}

\begin{theorem}[\cite{AAC18}]
\label{thdiamgenjohn}
If $n>k>i$ are non-negative integers, $n \geq 2k$ and $(n,k,i)\neq (2k,k,0)$ then 
\[
\diam(J(n,k,i))=\begin{cases}
    \displaystyle\left\lceil \frac{k-i-1}{n-2k+2i} \right\rceil + 1 
    &\text{if $n<3(k-i)-1$ or $i=0$;}\\
    3 &\text{if $3(k-i)-1  \leq n <3k-2i $ and $i \neq 0$;}\\
    \displaystyle\left\lceil \frac{k}{k-i} \right\rceil 
    &\text{if $n \geq 3k-2i  $ and $i \neq 0$.}
\end{cases}
\]

\end{theorem}

 Notice that given two vertices $A$ and $B$ in $J(n, k, i)$, they are adjacent if $|A \cap B|=i$, which means that $| \bar{A} \cap \bar{B}|=n-2k+i$. Thus, in such a case, $\bar{A}$ and $\bar{B}$ are adjacent in the graph $J(n, n-k, n-2k+i)$. Thus, the original result in  Theorem \ref{thdiamgenjohn} can be written for $n < 2k$ as follows.

\begin{remark}
\label{obsdiamj}
 If $n>k>i$ are non-negative integers, $n < 2k$ and $(n,k,i)\neq (2k,k,0)$ then
\[
\diam(J(n,k,i)) = \begin{cases}
    \makebox[9.6cm][s]{$\displaystyle \left\lceil \frac{k-i-1}{n-2k+2i} \right\rceil + 1$ \hfill if $n<3(k-i)-1$ or $n-2k+i=0$;} \\[1.5ex]
    \makebox[9.6cm][s]{$\displaystyle 3$ \hfill if $3(k-i)-1 \leq n <n+k-2i$ and $n \neq 2k-i$;} \\[1.5ex]
    \makebox[9.6cm][s]{$\displaystyle \left\lceil \frac{n-k}{k-i} \right\rceil$ \hfill if $2i \geq k$ and $n-2k+i \neq 0$.}
\end{cases}
\]
\end{remark}

In this paper, we are interested in {\em exact distance-$d$ Kneser graphs}, $K_{=d}(2k+r,k)$, with $k$ and $r$ positive integers.  In \cite{balak} the following theorem gives a complete characterization of the adjacency relation on $K_{=d}(2k+r,k)$ in terms of the cardinality of the intersection of its vertices, thus showing that they are another type of generalization of Kneser graphs. 

\begin{theorem}[\cite{balak}]
\label{thadj}
Let $k$ and $r$ be positive integers with $k\geq2$ and $1 \leq r < k-1$. Let $D= \diam(K(2k+r,k))$ and let $1 < d \leq D$ be a fixed integer. 
Two vertices $A$ and $B$ in $K_{=d}(2k+r,k)$ are adjacent if and only if

\begin{itemize}
\item[(i)] Case $d = 2p$ (with $p\geq 1$).  
\begin{enumerate}
\item $k-rp \leq |A \cap B| \leq k-rp + r -1$ when $d < D$. 
\item $rp - r +1 \leq |A \cap B| \leq k-rp+r-1$ if $d = D$. 
\end{enumerate}
\item[(ii)] Case $d = 2p+1$ (with $p\geq 1$). 
\begin{enumerate}
\item $rp-r+1 \leq |A \cap B| \leq rp$ when $d < D$.
\item $rp - r +1 \leq |A \cap B| \leq k-rp-1$ if $d = D$.
\end{enumerate}
\end{itemize}
\end{theorem}
 As in the different discussions of this article, we talk both about the distance between two vertices in $K(2k+r,k)$ and in $K_{=d}(2k+r,k)$, in order to avoid confusion and burdensome notation, we prefer to denote $\dist_K(A,B)$ instead of $\dist_{K(2k+r,k)}(A,B)$, and $\dist_{(K,d)}(A,B)$ instead of $\dist_{K_{=d}(2k+r,k)}(A,B)$, when the parameters of the Kneser graph and the exact distance-$d$ Kneser graph are clear from the context.
 In Section \ref{Sec: Distance}, we compute the distance of any two vertices in $K_{=d}(2k+r,k)$, and we prove the two following theorems.
\begin{theorem}
\label{thdistpar}
Let $A$ and $B$ be two non-adjacent vertices in $K_{=2p}(2k+r,k)$. If $|A\cap B|=s$ then
\[
\dist_{(K,2p)}(A,B)=\max\left\lbrace
2,
\left\lceil \frac{k-s}{rp}\right\rceil \right\rbrace.
\]
\end{theorem}
\begin{theorem}
\label{thdistimpar}
Let $A$ and $B$ be two non-adjacent vertices in $K_{=2p+1}(2k+r,k)$. 
If $|A\cap B|=s$ then
\[
\dist_{(K,2p+1)}(A,B)=
\min\left\lbrace 1+2\left\lceil\frac{\left|s-rp\right|}{2rp+r}\right\rceil , 2\left\lceil\frac{k-s}{2rp+r}\right\rceil\right\rbrace .
\]
\end{theorem}

Finally, in Section \ref{Sec: diameter}, we compute the diameter of graph $K_{=d}(2k+r,k)$, and prove the following theorem.
\begin{theorem}
    \label{thdiam}
Let $k$ and $r$ be positive integers with $k\geq2$ and $1 \leq r < k-1$. Let $D$ be the diameter of the Kneser graph $K(2k+r,k)$ and let $1 < d \leq D$. The diameter of the exact distance-$d$ Kneser graph $K_{=d}(2k+r,k)$ is given by
\[ 
\diam(K_{=d}(2k+r,k))=\begin{cases}
2 
& \text{$(i)$ if } d=D; \\
& \\
\left\lceil \frac{k}{rp}\right\rceil 
& \text{$(ii)$ if } d=2p<D ;\\
& \\
\left\lceil \frac{k+(r/2)}{2rp+r}\right\rceil 
& \text{$(iii)$ if } d=2p+1<D, \ 2\mid k-\frac{r}{2};\\
& \\
1+\left\lceil \frac{k+(r/2)-1}{2rp+r}\right\rceil 
& \text{$(iv)$ if } d=2p+1<D, \ 2\nmid k-\frac{r}{2}.
\end{cases}
\]
\end{theorem}


\section{Adjacency characterization}\label{Sec:Adjacency}

This section aims to provide some comments on the adjacency relation of vertices in $K_{=d}(2k+r,k)$, which will be useful in the next section.

For the remainder of this article, we denote $|A \cap B| = s$, $\diam(K(2k+r,k))=D$ and assume $d$ to be a fixed integer with $1 \leq d \leq D$.

Notice first that, by Proposition \ref{propvv},  $D=2$ when $r\geq k-1$ and thus, in this case, there are only two different distances: $1$ and $2$. As $K_{=1}(2k+r,k)$ is isomorphic to $K(2k+r,k)$ and $K_{=2}(2k+r,k)$ is isomorphic to the complement of $K(2k+r,k)$, we have the following remark.

\begin{remark}
\label{remarque1}
Let $k$ and $r$ be positive integers with $r\geq k-1$. Two vertices $A$ and $B$ in $K_{=2}(2k+r,k)$ are adjacent if and only if $1 \leq s \leq k-1$.
\end{remark}

The following remark is used in several proofs.

\begin{remark}
\label{laubonita}
Let $k$ and $r$ be positive integers with $k\geq2$ and $1 \leq r < k-1$, and let $k-1 = xr + y$ with $x,y$ integers such that $x > 0$ and $0 \leq y < r$. By Theorem \ref{thvv}, $D= \left\lceil \frac{k-1}{r} \right\rceil + 1$. Clearly, if $r$ divides $k-1$ then $y = 0$ and $D = x+1$, otherwise $D = x+2$. 
\end{remark}

As a direct consequence of Theorem \ref{thadj} we have the following result for $r=1$.
\begin{corollary}
Let $k$ be a positive integer with $k\geq2$ and $1 < d \leq k$ be a fixed integer. Two vertices $A$ and $B$ in $K_{=d}(2k+1,k)$ are adjacent if and only if
\[
s =\begin{cases}
k- \frac{d}{2} &\text{if $d$ is even;}\\
\frac{d-1}{2}&\text{if $d$ is odd.}
\end{cases}
\]
\end{corollary}


\section{Computing the distance function}\label{Sec: Distance}
In this section, we compute the distance function of two non-adjacent vertices in the exact distance-$d$ Kneser graph $K_{=d}(2k+r,k)$, where $k,r$ and $d$ are positive integers, with $r<k-1$. We assume $d \geq 2$, as the exact distance-$1$ Kneser graph is exactly the original Kneser graph whose distance function was computed in \cite{VV05}. By Remark \ref{remarque1}, we obtain the following.

\begin{remark}
\label{observation1}
Let $k$ and $r$ be positive integers with $r\geq k-1$. If $A$ and $B$ are two non-adjacent vertices in $K_{=2}(2k+r,k)$ then $\dist_{(K,2)}(A,B) = 2$.
\end{remark}

Let $A$ and $B$ be two vertices in $K_{=d}(2k+r,k)$. 
We denote by $C$ the set $A \cap B$ and we denote by $Z$ the set $[2k+r] \setminus C$. Notice that $|C| = s$ and $|Z| = r+s$. Thus, in this section, we use the following notation for sets $A,B,C$ and $Z$: $A = \{c_1,\ldots,c_s, a_1,\ldots, a_{k-s}\}$, $B = \{c_1,\ldots,c_s,b_1,\ldots,b_{k-s}\}$, $C = \{c_1,\ldots,c_s\}$, and $Z = \{z_1,\ldots,z_{r+s}\}$.

\subsection{Case \textit{d=2p}}

The proof of Theorem \ref{thdistpar} follows from Lemmas \ref{ldist2p1}, \ref{ldist2p2}, \ref{lem:Length of a path and intersection d=2p}, \ref{ldist2p2final} and Corollary \ref{cor:dist2p} given below.

\begin{lemma}
\label{ldist2p1}
Let $2p=D$. 
If $A$ and $B$ are two vertices in $K_{=2p}(2k+r,k)$ which are not adjacent, then $\dist_{(K,2p)}(A,B) = 2$.
\end{lemma}
\begin{proof}
As vertices $A$ and $B$ are not adjacent then $\dist_{(K,2p)}(A,B) \geq 2$. Moreover, by Theorem \ref{thadj} (Case $(i).2$), $s \leq rp-r$ or $s \geq k-rp+r$. 
\begin{itemize}
\item[(i)] 
Let $0 \leq s \leq rp-r$. As $rp-r+1 \leq k$, $rp \geq r-1$ and $k-rp+r-1-s \geq 0$, we have the following vertex $X$ in $K_{=2p}(2k+r,k)$: 
$$X=\{c_1,\ldots,c_s,a_1,\ldots,a_{rp-r+1-s},b_1,\ldots, b_{k-rp+r-1-s}, z_1,\ldots,z_s\}.$$
Clearly, $|A \cap X| = rp-r+1$ and $|B \cap X| = k-rp+r-1$. Thus, by Theorem \ref{thadj} (Case $(i).2$), $X$ is adjacent to $A$ and to $B$ in $K_{=2p}(2k+r,k)$. Therefore, $\dist_{(K,2p)}(A,B) \leq 2$. 
\item[(ii)] 
Let $k-rp+r \leq s \leq k-1$. As $rp-r+1 \leq k-rp+r-1$ and $|Z| = s+r > s \geq k-rp+r \geq rp-r+2 > rp-r+1$, we have the following vertex $X$ in $K_{=2p}(2k+r,k)$: 
$$X=\{c_1,\ldots,c_{k-rp+r-1},z_1,\ldots,z_{rp-r+1}\}.$$
Clearly, $|A \cap X| = |B \cap X| = k-rp+r-1$. Thus, by Theorem \ref{thadj} (Case $(i).2$), $X$ is adjacent to $A$ and to $B$ in $K_{=2p}(2k+r,k)$.  Therefore, $\dist_{(K,2p)}(A,B) \leq 2$. \qedhere
\end{itemize}
\end{proof}

\begin{lemma}
\label{ldist2p2}
Let $2p < D$. 
Let $A$ and $B$ be two vertices in $K_{=2p}(2k+r,k)$ which are not adjacent. 
If $s \geq k-rp+r$ then $\dist_{(K,2p)}(A,B)=2$. Otherwise, if $0 \leq s < k-rp$ then $\dist_{(K,2p)}(A,B) \leq \lceil \frac{k-s}{rp} \rceil$. 
\end{lemma}
\begin{proof}
By Theorem \ref{thadj} (Case $(i).1$), two vertices $A$ and $B$
are adjacent if and only if $k-rp \leq s \leq k-rp + r -1$. 
If $s \geq k-rp+r$ then as in the proof of Lemma \ref{ldist2p1} (Case $(ii)$), $\dist_{(K,2p)}(A,B) = 2$. Thus, assume that $0 \leq s \leq k-rp-1$. Let $k-rp-s = trp + q$, with $t\geq 0$ and $0 \leq q < rp$. We split the proof in two cases, either $t=0$, or $t>0$.
\begin{itemize}
\item $t=0$. In this case, $k-rp-s < rp$, which implies $k-s < 2rp$ and therefore, $2rp+s-k > 0$. In order to prove that $2rp+s-k \leq r+s$, it suffices to prove that $k \geq 2rp-r$. As in Remark \ref{laubonita}, we set $k-1 = xr +y$ with $x\geq 1$ and $0 \leq y \leq r-1$. 
\begin{itemize}
\item Assume $r$ divides $k-1$. In this case $D= x+1$. So, if $x+1 \geq 2p$ then $k \geq (2p-1)r + 1 > 2pr -r$.
\item Assume $r$ does not divide $k-1$. In this case, $D=x+2$. By theorem hypothesis, $x+2 > 2p$ and so, $k \geq (2p-1)r + y + 1 > 2pr-r$ because $y+1 \geq 2$.
\end{itemize}
Now we have the following vertex $X$ in $K_{=2p}(2k+r,k)$: 
$$X = \{c_1,\ldots,c_s, a_1,\ldots,a_{k-rp-s},b_1,\ldots,b_{k-rp-s}, z_1,\ldots,z_{2rp+s-k}\}.$$
Clearly, $|X \cap A| = |X \cap B| = k-rp$ and thus, $X$ is adjacent to both $A$ and $B$.

\item $t > 0$. Let $A = X_0$ and for $1 \leq i \leq t$, consider the following vertices in $K_{=2p}(2k+r,k)$:
$$X_i = \{c_1,\ldots,c_s,a_1,\ldots,a_{k-irp-s},b_1,\ldots,b_{irp}\}.$$
Notice that $|X_i \cap X_{i+1}| = s+ (irp) + (k - (i+1)rp -s) = k-rp$ and thus $X_i$ is adjacent to $X_{i+1}$ for $0 \leq i < t$. Moreover, $X_i$ is not adjacent to $X_j$, for $j>i+1$ or $j< i-1$ because $|X_i \cap X_j| < k-rp$. Now, we consider two cases:
\begin{itemize}
\item $rp$ divides $k-rp-s$. 

In this case, $X_t = \{c_1,\ldots,c_s, a_1,\ldots,a_{rp}, b_1,\ldots,b_{k-rp-s}\}$. Thus, $|X_t \cap B| = k-rp$ which implies that $X_t$ and $B$ are adjacent. Hence, $\dist_{(K,2p)}(A,B) \leq t+1$.

\item $rp$ does not divide $k-rp-s$.

\noindent In this case, $X_t = \{c_1,\ldots,c_s, a_1,\ldots,a_{rp+q}, b_1,\ldots,b_{k-rp-s-q}\}$. Thus, $|X_t \cap B| = k-rp-q < k-rp$ because $q>0$, which implies that $X_t$ and $B$ are not adjacent. 

Let us prove that $k-s-trp+rp-q \leq k-s$. As $t,q > 0$ then $rp \leq trp+q$ and so $k-s+rp - (trp+q) \leq k-s$.
Now we have the following vertex:
\[X_{t+1}=C \cup \{a_1,\ldots, a_q,a_{k-s-trp+1},\ldots, a_{k-s-trp+rp-q},\\b_1,\ldots,b_{trp+q}\}.\]
Notice that $|X_t \cap X_{t+1}| = s+trp +q = k-rp = |X_{t+1} \cap B|$ and thus, $X_{t+1}$ is adjacent to both $X_t$ and $B$. Finally, for this case, we obtain $\dist_{(K,2p)}(A,B) \leq t+2$.
\end{itemize}
\end{itemize}
Therefore, $\dist_{(K,2p)}(A,B) \leq \lceil \frac{k-rp-s}{rp} \rceil + 1 = \lceil \frac{k-s}{rp} \rceil$.
\end{proof}

In order to obtain a lower bound for the distance function when $d$ is even and $d < D$, we need the following results.

\begin{lemma}\label{lem:Length of a path and intersection d=2p}
Let $2p < D$. 
Let $A$ and $B$ be two vertices in $K_{=2p}(2k+r,k)$. If there is a path of length $\ell$ from $A$ to $B$, then $s\geq k-\ell rp$. 
\end{lemma}
\begin{proof}
    The proof follows by induction. 
    For the base case, if there is a path of length $1$, then $s\geq k-rp$ by Theorem \ref{thadj}.

    Assume there is a path of length $\ell$ from $A$ to $B$, and let $X$ be the vertex adjacent to $B$ in such path.
    As there is a path of length $(\ell-1)$ from $A$ to $X$,  by the inductive hypothesis, we have $|A\cap X|\geq k-(\ell-1)rp.$
    On the other hand, as $B$ and $X$ are adjacent, we have $|X\cap B|\geq k-rp$. 
    Then 
    \begin{align*}
        s=|A\cap B|\geq& |A\cap B\cap X|\\
        \geq &|A\cap X|-|X\cap \overline{B}|\\
        = &|A\cap X|-|X|+|X\cap B|\\
        \geq &
        k-(\ell-1)rp-k+k-rp\\
        =&k-\ell rp.
    \end{align*}
    
    Thus, the result follows by induction.
\end{proof}
\begin{corollary}\label{cor:dist2p}
Let $2p < D$. 
If $A$ and $B$ are two vertices in $K_{=2p}(2k+r,k)$ then
    $\dist_{(K,2p)}(A,B)\geq \frac{k-s}{rp}$.
\end{corollary}
\begin{proof}
There is a path of length $\dist_{(K,2p)}(A,B)$ from $A$ to $B$. Thus, applying Lemma \ref{lem:Length of a path and intersection d=2p} with $\ell=\dist_{(K,2p)}(A,B)$ yields, 
\[
    |A\cap B|\geq k-\dist_{(K,2p)}(A,B)rp.
\]
Solving for $\dist_{(K,2p)}(A,B)$ we get
\begin{align*}
    \dist_{(K,2p)}(A,B)\geq& \frac{k-|A\cap B|}{rp}.
\end{align*}
Thus, the result follows.
\end{proof}

\begin{lemma}
\label{ldist2p2final}
Let $2p < D$. 
Let $A$ and $B$ be two vertices in $K_{=2p}(2k+r,k)$ which are not adjacent. 
If $0 \leq s < k-rp$ then $\dist_{(K,2p)}(A,B) =\lceil \frac{k-s}{rp} \rceil$.
\end{lemma}
\begin{proof}
    The result follows from Lemma \ref{ldist2p2} and Corollary \ref{cor:dist2p}.
\end{proof}
To complete the proof of Theorem \ref{thdistpar}, we must show that $\left\lceil \frac{k-|A \cap B|}{rp}\right\rceil\leq 2$ whenever 
$
\dist_{(K,2p)}(A,B)=2$.
This is split into two cases. Namely, the case when $|A\cap B|\geq k-rp+r$ and the case when $d=D$. First, when $|A\cap B|\geq k-rp+r$, we have
\begin{align*}
    \frac{k-|A \cap B|}{rp}\leq &\frac{rp-r}{rp}
    \leq1.
\end{align*}
Thus, $\left\lceil \frac{k-|A \cap B|}{rp}\right\rceil\leq 1<2$.

When $d=D$ we have 
\[
p=\frac{D}{2}=\frac{1}{2}\left( \left\lceil\frac{k-1}{r}\right\rceil+1\right).
\]
Thus,
\begin{align*}
    \frac{k-|A \cap B|}{rp}=&2\left(\frac{k-|A \cap B|}{r\left\lceil\frac{k-1}{r}\right\rceil+r}\right)\\
    \leq& 2\left(\frac{k}{r\left(\frac{k-1}{r}\right)+r}\right)\\
    \leq& 2\left(\frac{k}{k}\right)\\
    =&2.
\end{align*}
Similarly, this implies $\left\lceil \frac{k-|A \cap B|}{rp}\right\rceil\leq 2$.

\subsection{Case \textit{d = 2p+1}}

The proof of Theorem \ref{thdistimpar} follows from Lemmas \ref{ldist-imp1}, \ref{ldist-imp2}, \ref{lem:Length of a path and intersection} and \ref{l:distimpcasolargo}, and Corollary \ref{cor dist 2 o 3 impar} given below.

\begin{lemma}
\label{ldist-imp1}
Let  $2p+1=D$ and $r$ does not divide $k-1$. 
If $A$ and $B$ are two vertices in $K_{=2p+1}(2k+r,k)$ which are not adjacent, then $\dist_{(K,2p+1)}(A,B) = 2$.
\end{lemma}
\begin{proof}
As vertices $A$ and $B$ are not adjacent, $\dist_{(K,2p+1)}(A,B) \geq 2$. Moreover, by Theorem \ref{thadj} (Case $(ii).2$), $|A \cap B| \leq rp-r$ or $|A \cap B| \geq k-rp$. As in Remark \ref{laubonita}, we set $k-1 = xr + y$, with $0 < y < r$. Since $r$ does not divide $k-1$, it follows that  $D = x+2$.
So, in this case, $x+2 =2p+1$ which implies that $x = 2p-1$ and therefore, $k = 2rp - r + y+1$.  
\begin{itemize}
\item[(i)] Case $0 \leq |A \cap B| \leq rp-r$. As $r,y>0$, $s \leq rp - r$, and $k-rp-1-s \geq rp-r+y-s \geq y$, the quantities $rp - r + 1 -s$ and $k - rp -1 - s$ are both positive. Thus, we can consider the vertex $X$ in $K_{=2p+1}(2k+r,k)$: 
$$X=C \cup \{a_1,\ldots,a_{rp-r+1-s},b_1,\ldots, b_{k-rp-1-s}, z_1,\ldots,z_{s+r}\} .$$
Clearly, $|A \cap X| = rp-r+1$ and $|B \cap X| = k-rp-1$. Hence, by Theorem \ref{thadj} (Case $(ii).2$), $X$ is adjacent to $A$ and to $B$ in $K_{=2p+1}(2k+r,k)$.  Therefore, $\dist_{(K,2p+1)}(A,B) \leq 2$. 
\item[(ii)] Case $k-rp \leq |A \cap B| \leq k-1$. As $k = 2rp -r +y +1$, we have $0 < k-rp-1 < s$ and $s+r > rp +1$. Thus, we can consider the following vertex $X$ in $K_{=2p+1}(2k+r,k)$: 
$$X=\{c_1,\ldots,c_{k-rp-1},z_1,\ldots,z_{rp+1}\}.$$
Clearly, $|A \cap X| = |B \cap X| = k-rp-1$. Hence, by Theorem \ref{thadj} (Case $(ii).2$), $X$ is adjacent to $A$ and to $B$ in $K_{=2p+1}(2k+r,k)$.  Therefore, $\dist_{(K,2p+1)}(A,B) \leq 2$. \qedhere
\end{itemize}
\end{proof}

\begin{lemma}
\label{ldist-imp2}
Let $2p+1=D$ and $r$ divides $k-1$.
If $A$ and $B$ are two vertices in $K_{=2p+1}(2k+r,k)$ which are not adjacent, then $\dist_{(K,2p+1)}(A,B) = 2$.
\end{lemma}
\begin{proof}
As vertices $A$ and $B$ are not adjacent then $\dist_{(K,2p+1)}(A,B) \geq 2$. Moreover, by Theorem \ref{thadj} (Case $(ii).1$), $s \leq rp-r$ or $s \geq rp+1$. Recall that, in this case, using Remark \ref{laubonita}, we obtain $D = x+1$ and thus $k=2rp + 1$.
\begin{itemize}
\item[(i)] Case $0 \leq s \leq rp-r$. This case is similar to the Case $(i)$ in Lemma \ref{ldist-imp1}.
\item[(ii)] Case $rp+1 \leq s \leq k-1$. Notice that $k-rp = rp+1<s+r$. Thus, we can consider the following vertex $X$ in $K_{=2p+1}(2k+r,k)$: 
$$X=\{c_1,\ldots,c_{rp},z_1,\ldots,z_{k-rp}\}.$$
Clearly, $|A \cap X| = |B \cap X| = rp$. Hence, by Theorem \ref{thadj} (Case $(ii).1$), $X$ is adjacent to $A$ and to $B$ in $K_{=2p+1}(2k+r,k)$. Therefore, $\dist_{(K,2p+1)}(A,B) \leq 2$. \qedhere
\end{itemize}
\end{proof}

\begin{lemma}
\label{ldist-imp3}
Let $2p+1<D$ and let $A$ and $B$ be two vertices in $K_{=2p+1}(2k+r,k)$. If $s \leq rp-r$, then
\[
\dist_{(K,2p+1)}(A,B)=\begin{cases}
    2 &\text{if $s\geq k-2rp-r$;}\\
    3 &\text{if $s< k-2rp-r$.}
\end{cases}
\]
\end{lemma}
\begin{proof}
\begin{itemize}
\item[(i)] Case $s\geq k-2rp -r$. As $s\leq rp-r$, by Theorem \ref{thadj},  $A$ and $B$ are not adjacent. Thus, for this case, it is enough to show that $\dist_{(K,2p+1)}(A,B)\leq 2$. Notice that the set 
\[ X=\{a_1,\ldots,a_{rp},b_1,\ldots,b_{rp},z_1,\ldots,z_{k-2rp}\}\]
is a vertex in $K_{=2p+1}(2k+r,k)$ because $k > 2rp$ when $d<D$.
Clearly, $|A \cap X| = |X \cap B| = rp$ and, by Theorem \ref{thadj} (Case $(ii).1$), $X$ is adjacent to $A$ and to $B$ in $K_{=2p+1}(2k+r+k)$. Therefore, 
$\dist_{(K,2p+1)}(A,B) = 2$. 
\item[(ii)] Case $s<k-2rp -r$.
The proof of this case is split into two parts. First we show that $\dist_{(K,2p+1)}(A,B)\leq 3$. 
As  $rp -r +1 - s > 0$ because $s \leq rp-r$ we can consider the following vertices $X_1,X_2$ in $K_{=2p+1}(2k+r,k)$: 
$$X_1=\{c_1,\ldots,c_s,a_1,\ldots,a_{rp-r+1-s},b_1,\ldots,b_{k-rp+r-1}\},$$
$$X_2=\{c_1,\ldots,c_s,a_{rp-r+1-s+1},\ldots,a_{k-s},b_1,\ldots,b_{rp-r+1-s}\}.$$

Notice that  $|X_1\cap X_2|=|A \cap X_1| = |X_2 \cap B| = rp-r+1$. Thus, by Theorem \ref{thadj}, $X_1$ is adjacent to $A$ and to $X_2$, and $B$ is adjacent to $X_2$. Hence,  $\dist_{(K,2p+1)}(A,B) \leq 3$. 

To complete the proof, we need to show that no vertex $Y$ can be adjacent to both $A$ and $B$. As $|A\cap B|< k-r-2rp$, we have $|Z|=r+s< k-2rp$. Thus, if $|Y|=k$, then $|Y\cap (A\cup B)|\geq |Y|-|Z| \geq 2rp+1$. This implies that $\max\{|Y\cap A|,|Y\cap B|\}\geq rp+1$. Hence, by Theorem \ref{thadj}, $Y$ cannot be simultaneously adjacent to $A$ and $B$. Therefore, $\dist_{(K,2p+1)}(A,B)\geq 3$.\qedhere
\end{itemize}
\end{proof}

\begin{corollary}\label{cor dist 2 o 3 impar}
Let $A$ and $B$ be two vertices in $K_{=2p+1}(2k+r,k)$. If $2p+1=D$ or $s\leq rp-r$, then
\[
    \dist_{(K,2p+1)}(A,B)=\min\left\lbrace 1+2\left\lceil\frac{\left|s-rp\right|}{2rp+r}\right\rceil , 2\left\lceil\frac{k-s}{2rp+r}\right\rceil\right\rbrace.
\]
\end{corollary}
\begin{proof}
Let $m=\min\left\lbrace 1+2\left\lceil\frac{\left|s-rp\right|}{2rp+r}\right\rceil , 2\left\lceil\frac{k-s}{2rp+r}\right\rceil\right\rbrace$.
By Lemma \ref{ldist-imp3}, we need to show that
\[
m =\begin{cases}
2&\text{$(i)$ if } d=D,\\
& \text{$(ii)$ if $d< D$, } k-2rp-r\leq s\leq rp-r;\\
3&\text{if $d< D$, $s\leq rp-r$, and $s< k-2rp-r$.}
\end{cases}
\]
Notice that  $s \not\in \{rp,k\}$, as $A$ and $B$ are not adjacent and $A\neq B$. This implies that 
\[
1+2\left\lceil\frac{\left|s-rp\right|}{2rp+r}\right\rceil\geq 3,
\]
and that
\[
2\left\lceil\frac{k-s}{2rp+r}\right\rceil\geq 2.
\]
Thus, to prove that 
$m=2$, it suffices to show that $\left\lceil\frac{k-s}{2rp+r}\right\rceil\leq 1$.

When $d=D$, we have $2p+1=\left\lceil\frac{k-1}{r}\right\rceil+1$. Then
\[
p=\frac{1}{2}\left\lceil\frac{k-1}{r}\right\rceil\geq \frac{k-1}{2r}.
\]
Thus, 
\begin{align*}
    \left\lceil\frac{k-s}{2rp+r}\right\rceil \leq & \left\lceil\frac{k}{k-1+r}\right\rceil
    =  1.
\end{align*}
When $d<D$ and $k-2rp-r\leq s\leq rp-r$, 
\begin{align*}
\left\lceil\frac{k-s}{2rp+r}\right\rceil 
\leq&  \left\lceil\frac{k-k+2rp+r}{2rp+r}\right\rceil 
=  1.
\end{align*}

Thus, in these cases $m=2$.

Assume now that $d<D$, $s<k-2rp-r$ and $s\leq rp-r$.
Notice that
\begin{align*}
\frac{k-s}{2rp+r} 
\geq&  \frac{k-k+2rp+r+1}{2rp+r} \\
=&  1+\frac{1}{2rp+r}.
\end{align*}
Thus,
\begin{align*}
2\left\lceil\frac{k-s}{2rp+r} \right\rceil
\geq&  2\left\lceil 1+\frac{1}{2rp+r}\right\rceil
= 4.
\end{align*}
But, as $s\leq rp-r$, 
\begin{align*}
    1+2\left\lceil\frac{\left|s-rp\right|}{2rp+r}\right\rceil=& 1+2\left\lceil\frac{rp-s}{2rp+r}\right\rceil\\
    \leq & 1 +2\left\lceil\frac{rp}{2rp+r}\right\rceil\\
    \leq &3.
\end{align*}
Therefore, in this case, $m=3$.
\end{proof}

For the remaining cases, we first use the existence of a path of a certain length between $A$ and $B$ to give a lower bound for $|A\cap B|$. Then, we use this bound to find the actual distance.
\begin{lemma}\label{lem:Length of a path and intersection}
    Let $2p+1<D$. 
    If there is a path of length $2\ell$ from $A$ to $B$, then $|A\cap B|\geq k-2\ell rp-r\ell$. If there is a path of length $2\ell +1$ from $A$ to $B$, then $|A\cap B|\leq (2\ell+1)rp+r\ell$.
\end{lemma}
\begin{proof}
The proof follows by induction. 
    For the base case, when $\ell=0$, if there is a path of length $1=2\cdot 0+1$, then $|A\cap B|\leq (2\ell+1)rp+r\ell$ by Theorem \ref{thadj}.

Assume there is a path of length $2\ell$ from $A$ to $B$, and let $X$ be the vertex adjacent to $B$ in such a path. As there is a path of length $2 (\ell-1)+1$ from $A$ to $X$,  by the inductive hypothesis, we have 
\begin{align*}        
    |A\cap X|\leq& [2(\ell-1)+1]rp+r(\ell-1)\\
    =&(2\ell -1)rp+r(\ell-1).
\end{align*}
Thus, we have an upper bound for the size of the complement of $A\cup X$, 
\begin{align*}
        |\overline{A\cup X}|= & 2k+r -|A|-|X|+|A\cap X|\\
        \leq& 2k+r-k-k+(2\ell -1)rp+r(\ell-1)\\
        =&(2\ell -1)rp+r\ell.
\end{align*}
On the other hand, as $B$ and $X$ are adjacent, by Theorem \ref{thadj} we have 
\[
    rp-r+1\leq |B\cap X|\leq rp.
\]
Hence,
\begin{align*}
        |B\cap \overline{A}|=&|B\cap (\overline{A\cup X})|+ |B\cap \overline{A}\cap X|\\
        \leq& |\overline{A\cup X}|+ |B\cap X|\\
        \leq & (2\ell -1)rp+r\ell+rp\\
        =& 2\ell rp+r\ell.
\end{align*}
Thus, 
\begin{align*}
        |B\cap A|=&k-|B\cap \overline{A}|\\
        \geq & k-2\ell rp-r\ell.\\
\end{align*}
Which proves the even case.

For the odd case, assume there is a path of length $2\ell+1$ from $A$ to $B$, and let $X$ be the vertex adjacent to $B$ in such a path. As there is a path of length $2 \ell$ from $A$ to $X$,  by the inductive hypothesis, we have 
\begin{align*}        
    |A\cap X|\geq& k-2\ell rp-r\ell.
\end{align*}
Thus, we have an upper bound for the size of $A\cap \overline{X}$,
\begin{align*}
        |A\cap \overline{X}|\leq& k-(k-2\ell rp-r\ell)\\
        =&2\ell rp+r\ell.
\end{align*}
On the other hand, as $B$ and $X$ are adjacent, by Theorem \ref{thadj} we have 
\[
    rp-r+1\leq |B\cap X|\leq rp.
\]
Hence
\begin{align*}
        |A\cap B|=&|A\cap B\cap \overline{X}|+ |A\cap B\cap X|\\
        \leq& |A\cap \overline{X}|+ |B\cap X|\\
        \leq & 2\ell rp+r\ell + rp\\
        =& (2\ell+1) rp+r\ell.
\end{align*}

Therefore, the result follows by induction.
\end{proof}
By Lemma \ref{lem:Length of a path and intersection}, if $s\geq rp+1$, then 
$\dist_{(K,2p+1)}(A,B)$ is given by
\[
\min \left(\{2\ell+1 \ / \  (2\ell+1) rp+r\ell\geq s\}\cup \{2\ell \ / \ k-2\ell rp-r\ell\leq s\}\right).
\]
Assume $d=\min \{\ell \ / \  (2\ell+1) rp+r\ell\geq s\}$. Then 
\begin{align*}
    (2(d-1)+1) rp+r(d-1)< s\leq (2d+1) rp+rd.
\end{align*}
Which yields
\begin{align*}
    \frac{s-rp}{2rp+r}\leq d <\frac{s+rp+r}{2rp+r}.
\end{align*}
But, as $d$ is an integer, 
\begin{align*}
    \left\lceil\frac{s-rp}{2rp+r}\right\rceil\leq d &<\left\lceil\frac{s+rp+r}{2rp+r}\right\rceil\\
        \left\lceil\frac{s-rp}{2rp+r}\right\rceil\leq d &<\left\lceil\frac{s-rp+2rp+r}{2rp+r}\right\rceil\\
                \left\lceil\frac{s-rp}{2rp+r}\right\rceil\leq d &<\left\lceil\frac{s-rp}{2rp+r}+1\right\rceil
                \\
                \left\lceil\frac{s-rp}{2rp+r}\right\rceil\leq d &<\left\lceil\frac{s-rp}{2rp+r}\right\rceil+1.
\end{align*}
Which implies,  $d=\left\lceil\frac{s-rp}{2rp+r}\right\rceil$, and $2d+1=2\left\lceil\frac{s-rp}{2rp+r}\right\rceil +1$.

Similarly, if $d=\min\{\ell \ / \ k-2\ell rp-r\ell\leq s\}$,  we have
\begin{align*}
    k-2drp -rd\leq s<k-2(d-1)rp -r(d-1),
\end{align*}
which, solving for $d$, implies
\begin{align*}
    \frac{k-s}{2rp+r}\leq d<\frac{k-s}{2rp+r} +1.
\end{align*}
Thus, $d=\left\lceil\frac{k-s}{2rp+r}\right\rceil$, and $2d=2\left\lceil\frac{k-s}{2rp+r}\right\rceil$.
Hence, the distance between $A$ and $B$ is
\[
\dist_{(K,2p+1)}(A,B) =\min\left\lbrace 1+2\left\lceil\frac{s-rp}{2rp+r}\right\rceil , 2\left\lceil\frac{k-s}{2rp+r}\right\rceil\right\rbrace.
\]
We have proved the following.
\begin{lemma}\label{l:distimpcasolargo}
    Let $2p+1<D$ and let $A$ and $B$ be two vertices in $K_{=2p+1}(2k+r,k)$. If $s\geq rp+1$,
    then $\dist_{(K,2p+1)}(A,B) =\min\left\lbrace 1+2\left\lceil\frac{s-rp}{2rp+r}\right\rceil , 2\left\lceil\frac{k-s}{2rp+r}\right\rceil\right\rbrace.$
\end{lemma}
\begin{proof}
    The result follows from the preceding discussion.
\end{proof}


\section{Computing the diameter}\label{Sec: diameter}
Notice first that, by Remark \ref{observation1}, we have the following result.

\begin{remark}
\label{observation2}
If $k$ and $r$ are positive integers with $k\geq2$ and $r \geq k-1$, then, $\diam(K_{=2}(2k+r,k))=2$.
\end{remark}
We can now proceed to the proof of Theorem \ref{thdiam}.


\begin{proof}[Proof of Theorem \ref{thdiam}]
    Case $(i)$ follows from Lemmas \ref{ldist2p1}, \ref{ldist-imp1}, and \ref{ldist-imp2}. 
    
    Case $(ii)$ follows from Lemmas \ref{ldist2p2} and \ref{ldist2p2final} by having $A\cap B=\emptyset$. 
    
    For case $(iii)$, assume $d=2p+1<D$. 
    By Lemma \ref{l:distimpcasolargo}
    we have 
    \[
\diam(K_{=d}(2k+r,k))=\max_{s}\min \left\lbrace 1+2\left\lceil\frac{s-rp}{2rp+r}\right\rceil , 2\left\lceil\frac{k-s}{2rp+r}\right\rceil\right\rbrace .
    \]
    Let $A$ and $B$ be two vertices in $K_{=d}(2k+r,k)$, and notice
    \[
    \dist_{(K,d)}(A,B)\leq 1+2\left\lceil\frac{s-rp}{2rp+r}\right\rceil < 3+2\left(\frac{s-rp}{2rp+r}\right)
    \]
    and 
    \[
    \dist_{(K,d)}(A,B)\leq 2\left\lceil\frac{k-s}{2rp+r}\right\rceil < 2+2\left(\frac{k-s}{2rp+r}\right).
    \] 
    We find the value of $s$ that maximizes the distance by doing
    \begin{align*}
        3+2\left(\frac{s-rp}{2rp+r}\right)=&2+2\left(\frac{k-s}{2rp+r}\right)\\
        2s=&k+rp-\frac{2rp+r}{2}\\
        s=&\frac{k-(r/2)}{2}.   
    \end{align*}
    In the case when $\frac{k-(r/2)}{2}$ is an integer, it is a valid value to take for $s$. In such a case $\diam(K_{=d}(2k+r,k))$ is given by
    \begin{align*}
& \min\left\lbrace 1+2\left\lceil\frac{(k-(r/2))/2-rp}{2rp+r}\right\rceil , 2\left\lceil\frac{k-(k-(r/2))/2}{2rp+r}\right\rceil\right\rbrace .\\ 
\end{align*}  
Letting $k=(4rp+2r)a+b$, with $0\leq b<4rp+2r$
\begin{align*}
2\left\lceil\frac{k-(k-(r/2))/2}{2rp+r}\right\rceil=&2\left\lceil\frac{(k+(r/2))/2}{2rp+r}\right\rceil\\
    =&2a+2\left\lceil\frac{b+(r/2)}{4rp+2r}\right\rceil \\
    =&\begin{cases}
        2a+2&\text{if }b\leq 4rp+r+(r/2)\\
        2a+4&\text{if }b>4rp+r+(r/2),\\
    \end{cases}
\end{align*}
and
\begin{align*}
1+2\left\lceil\frac{(k-(r/2))/2-rp}{2rp+r}\right\rceil &=1+2\left\lceil\frac{k-(r/2)-2rp}{4rp+2r}\right\rceil\\
    &=1+2(a-1)+2\left\lceil\frac{4rp+2r+b-(r/2)-2rp}{4rp+2r}\right\rceil\\
    &=1+2(a-1)+2\left\lceil\frac{2rp+b+(3r/2)}{4rp+2r}\right\rceil\\
    =&\begin{cases}
        2a+1&\text{if }b\leq 2rp+(r/2)\\
        2a+3&\text{if }b>2rp+(r/2).\\
    \end{cases}
\end{align*}
Thus, when $\frac{k-(r/2)}{2}$ is an integer, we have
    \begin{align*}
\diam(K_{=d}(2k+r,k))=2a+\begin{cases}
1&\text{if }b\leq 2rp+\lfloor(r/2)\rfloor\\
2&\text{if } 2rp+\lfloor(r/2)\rfloor<b\leq 4rp+r+\lfloor (r/2)\rfloor \\
3&\text{if } 4rp+r+\lfloor (r/2)\rfloor<b. \\
\end{cases}
\end{align*}
Notice that in this case 
\begin{align*}
    \left\lceil \frac{k+(r/2)}{2rp+r}\right\rceil=&\left\lceil \frac{2a(2rp+r)+b+(r/2)}{2rp+r}\right\rceil \\
    =&2a+\left\lceil \frac{b+(r/2)}{2rp+r}\right\rceil \\
    =&2a+\begin{cases}
        1&\text{if } b\leq 2rp +\lfloor (r/2)\rfloor \\
        2&\text{if } 2rp +\lfloor (r/2)\rfloor <b\leq 4rp +r+\lfloor (r/2)\rfloor \\
        3&\text{if }4rp +r+\lfloor (r/2)\rfloor<b.
    \end{cases}
\end{align*}
Thus, when $2p+1=d<D$ and $\frac{k-(r/2)}{2}$ is an integer, 
\[\diam(K_{=d}(2k+r,k))=\left\lceil \frac{k+(r/2)}{2rp+r}\right\rceil.\]

   For case $(iv)$, when $\frac{k-(r/2)}{2}$ is not an integer, notice that 
    $1+2\left\lceil\frac{s-rp}{2rp+r}\right\rceil$ and $2\left\lceil\frac{k-s}{2rp+r}\right\rceil$ are a non-decreasing and a non-increasing function of $s$, respectively. Thus, for $1+2\left\lceil\frac{s-rp}{2rp+r}\right\rceil$ we 
    consider the value $s=
    \frac{k-(r/2)-1}{2}$ 
and for $2\left\lceil\frac{k-s}{2rp+r}\right\rceil$ the 
value  $s=\frac{k-(r/2)+1}{2}$. As the diameter is the maximum of the distances, in this case, $\diam(K_{=d}(2k+r,k))$ is given by
\begin{align*}
& \max\left\lbrace 1+2\left\lceil\frac{(k-(r/2)-1)/2-rp}{2rp+r}\right\rceil , 2\left\lceil\frac{k-(k-(r/2)+1)/2}{2rp+r}\right\rceil\right\rbrace . \\ 
\end{align*}
Letting $k=(4rp+2r)a+2b+\xi$, with $0\leq 2b+\xi <4rp+2r$ and $0\leq \xi\leq 1$ and $r=4q+2t+\zeta$, with $0\leq t,\zeta\leq 1$ and $q = \lfloor\frac{r}{4}\rfloor$, we have
\begin{align*}
&\left\lceil\frac{(k-(r/2)-1)/2-rp}{2rp+r}\right\rceil= 
\left\lceil\frac{
(2rp+r)a+b-q-rp+(\xi-t-(\zeta/2)-1)/2}{2rp+r}\right\rceil
\\
   &=(a-1)+\left\lceil\frac{b+2rp+4q+2t+\zeta-q-rp+(\xi-t-(\zeta/2)-1)/2}{2rp+r}\right\rceil\\
   &=(a-1)+\left\lceil\frac{b+rp+3q+2t+\zeta+(\xi-t-(\zeta/2)-1)/2}{2rp+r}\right\rceil.\\
   \end{align*}
   Thus,
\begin{align*}
1+2\left\lceil\dfrac{\dfrac{k-(r/2)-1}{2}-rp}{2rp+r}\right\rceil
&=1+2(a-1)+2\left\lceil\frac{b+rp+3q+2t+\zeta+\dfrac{\xi-t-(\zeta/2)-1}{2}}{2rp+r}\right\rceil\\
    &=\begin{cases}
        2a+1&\text{if }b\leq rp+q-(\xi-t-(\zeta/2)-1)/2\\
        2a+3&\text{if }b>rp+q-(\xi-t-(\zeta/2)-1)/2
    \end{cases}
\end{align*}
and
\begin{align*}
2\left\lceil\dfrac{k-\dfrac{k-(r/2)+1}{2}}{2rp+r}\right\rceil
=& 2\left\lceil\frac{(4rp+2r)a+2b+\xi-\dfrac{(4rp+2r)a+2b+\xi-2q-t-(\zeta/2)+1}{2}}{2rp+r}\right\rceil\\
=&2\left\lceil\frac{(2rp+r)a+b+\xi+q-(\xi-t-(\zeta/2)+1)/2}{2rp+r}\right\rceil\\
    =&2a+2\left\lceil\frac{b+q+(\xi+t+(\zeta/2)-1)/2}{2rp+r}\right\rceil\\
    =&\begin{cases}
        2a+2&\text{if }b\leq 2rp+3q+2t+\zeta-(\xi+t+\frac{\zeta}{2}-1)/2\\
        2a+4&\text{if }b>2rp+3q+2t+\zeta-(\xi+t+\frac{\zeta}{2}-1)/2.\\
    \end{cases}
\end{align*}
Thus,
$\diam(K_{=d}(2k+r,k)) =2a+\alpha$, where
\begin{align*}    
\alpha= \begin{cases}
2&\text{if }b\leq rp+q-\frac{\xi}{2}+\frac{t}{2}+\frac{\zeta}{4}+\frac{1}{2}\\
3&\text{if }rp+q-\frac{\xi}{2}+\frac{t}{2}+\frac{\zeta}{4}+\frac{1}{2}< b\leq 2rp+3q+2t-\frac{\xi}{2}-\frac{t}{2}+\frac{3\zeta}{4}+\frac{1}{2}\\  4&\text{if }b>2rp+3q+2t-\frac{\xi}{2}-\frac{t}{2}+\frac{3\zeta}{4}+\frac{1}{2}.\\
\end{cases}
\end{align*}
But notice that 
\begin{align*}
   &1+ \left\lceil \frac{k+(r/2)-1}{2rp+r}\right\rceil =1+\left\lceil \frac{(2rp+r)2a+2b+\xi+2q+t+(\zeta/2)-1}{2rp+r}\right\rceil\\
    &=1+2a+\left\lceil \frac{2b+\xi+2q+t+(\zeta/2)-1}{2rp+r}\right\rceil\\
    &=2a+\begin{cases}
        2&\text{if }2b\leq 2rp+2q+t+(\zeta/2)-\xi+1\\
3&\text{if }2rp+2q+t+\frac{\zeta}{2}-\xi+1< 2b\leq 4rp+6q+3t+3\frac{\zeta}{2}-\xi+1\\  
4&\text{if }2b>4rp+6q+3t+3(\zeta/2)-\xi+1\\
    \end{cases}\\
    &=2a+\begin{cases}
        2&\text{if }b\leq rp+q-\frac{\xi}{2}+\frac{t}{2}+\frac{\zeta}{4}+\frac{1}{2}\\
3&\text{if }rp+q-\frac{\xi}{2}+\frac{t}{2}+\frac{\zeta}{4}+\frac{1}{2}< b\\
&\text{and } b\leq 2rp+3q+2t-\frac{\xi}{2}-\frac{t}{2}+\frac{3\zeta}{4}+\frac{1}{2}\\  
4&\text{if }b>2rp+3q+2t-\frac{\xi}{2}-\frac{t}{2}+\frac{3\zeta}{4}+\frac{1}{2}.\\
    \end{cases}
\end{align*}
Thus, in this case,
\[
\diam(K_{=d}(2k+r,k))= 1+\left\lceil \frac{k+(r/2)-1}{2rp+r}\right\rceil .
\]

\end{proof}


\section{Concluding remarks}
Notice that Theorems \ref{thvv}, \ref{thdiamgenkn}, \ref{thdiamgenjohn}
and \ref{thdiam} give the diameters of ``Kneser-like'' graphs where two vertices are adjacent if the cardinality of the intersection of the corresponding sets is either a fixed value, upper bounded, or both upper and lower bounded, although the latter with some specific bounds. It is natural then to study the diameter when other restrictions are applied to said cardinality. As a general problem, we propose the following.
\begin{problem}\label{problem}
    Let $c$ and $C$ be non-negative integers with $c\leq C\leq k$. What is the diameter of a graph $G$ with vertex set $V(G)=[n]^k$ in each of the following cases?
\begin{enumerate}[(i)]
        \item Two vertices $A,B$ are adjacent if $c\leq |A\cap B|\leq C$.
        \item Two vertices $A,B$ are adjacent if $|A\cap B|\leq c$ or $|A\cap B|\geq C$.
    \end{enumerate}
\end{problem}

We note that several cases of Problem~\ref{problem}~(i) are addressed in the following results:
\begin{itemize}
    \item Theorem~\ref{thvv} covers the case $c = C = 0$.
    \item Theorem~\ref{thdiamgenkn} covers the case $c = 0$.
    \item Theorem~\ref{thdiamgenjohn} covers the case $c = C$.
    \item Theorem~\ref{thdiam} addresses additional, more general cases, though still under certain restrictions on $c$ and $C$. Specifically:
    \begin{itemize}
        \item Cases~$(i)$ 1. and $(ii)$ 1.: some instances where $c \equiv 1 \text{ or } k \pmod{n - 2k}$ and $C = c+ n - 2k - 1$;
        \item Cases~$(i)$ 2. and $(ii)$ 2.: some instances where $c \equiv 1 \pmod{n - 2k}$ and $C  \equiv k - 1 \pmod{n - 2k}$;
        \end{itemize}
    We emphasize that not all parameter choices satisfying the conditions above are necessarily covered by Theorem \ref{thdiam}.
\end{itemize}

Some other cases can be studied through power graphs. On one hand, using the power graph of the Johnson graph (the generalized Johnson graph with $i=k-1$) we can obtain an answer for Problem \ref{problem} $(i)$ when $C=k$. On the other hand, using the power graph of the Kneser graph we can obtain some results for Problem \ref{problem} $(ii)$.

Starting with generalized Johnson graphs, by Theorem \ref{thdistgenjohn}, $A$ and $B$ are adjacent in $J(n,k,k-1)^{k-l}$ if and only if $|A\cap B| \geq k-l$, whereas, by Lemma \ref{lvv},
they are adjacent in $K(n,k)^p$ if and only if 
$\frac{p(n-2k)}{2}\leq |A\cap B|$ or
$|A\cap B|\leq k-\frac{p(n-2k)}{2}$. Thus, they cover the case where adjacency corresponds to a lower bound, and when it corresponds to the union of an upper bounded set and a lower bounded set, although the latter with some specific bounds. Finding the diameter of a power graph, $\Gamma^p$, in terms of the diameter of the original graph, $\Gamma$, is a simple exercise, as
\[
\dist_{\Gamma^p}(A,B) = q \iff (q-1)p \leq \dist_{\Gamma}(A,B)\leq pq.
\] Thus, 
\[
\diam(\Gamma^p)=\left\lceil \frac{\diam(\Gamma)}{p} \right\rceil.
\]
Combining this with Theorems \ref{thvv} and \ref{thdiamgenjohn} yields the following.
\begin{theorem}\label{thdiampower}
    Let $k$, $r$ and $p$ be positive integers. Then
    \[\diam(K(2k+r,k)^p)= \left\lceil  \frac{\left\lceil \frac{k-1}{r} \right\rceil + 1}{p} \right\rceil\\
    \]
    and
    \[\diam(J(n,k,k-1)^{k-l})=  \left\lceil  \frac{n-k}{k-l} \right\rceil.
    \]
\end{theorem}

Although Theorems \ref{thdiam} and \ref{thdiampower} present the diameter of graphs when two vertices are adjacent if their intersection is both upper and lower bounded (the case of exact distance Kneser graphs) and when they are adjacent if the intersection is the union of an upper bounded and a lower bounded set (the power of Kneser graphs), these bounds are not indepedent of each other. Given any two constants $c<C$, it would be interesting to study Problem \ref{problem} when no restrictions are applied to $c$ and $C$.

 In terms of the Wiener index, we can fix a vertex $A$ in $K_{=d}(2k+r,k)$, and, 
using Theorems \ref{thdistpar} and \ref{thdistimpar},
study their distance to every vertex in the graph. These yield  the following sums
\begin{align*}
\begin{cases} \displaystyle\sum_{s=0}^{k-1}\binom{k}{s}\binom{n-k}{k-s}\max\left\lbrace
2,
\left\lceil \frac{k-s}{rp}\right\rceil \right\rbrace & \text{if $d=2p$}\\[15pt]
\displaystyle\sum_{s=0}^{k-1}\binom{k}{s}\binom{n-k}{k-s}\min\left\lbrace 1+2\left\lceil\frac{\left|s-rp\right|}{2rp+r}\right\rceil , 2\left\lceil\frac{k-s}{2rp+r}\right\rceil\right\rbrace & \text{if $d=2p+1$,}\\
    \end{cases}
\end{align*}
where $\binom{k}{s}\binom{n-k}{k-s}$ is counting the number of vertices that intersect $A$ in $s$ elements. 
Afterwards, we can multiply by the number of vertices in the graph, $\binom{n}{k}$ and divide by $2$, as each edge would be counted twice.
Thus, the following formulas for the Wiener index $W(K_{=d}(2k+r,k))$ can be obtained:
\begin{align*}
    \begin{cases}
\displaystyle\frac{1}{2}\binom{n}{k}\sum_{s=0}^{k-1}\binom{k}{s}\binom{n-k}{k-s}\max\left\lbrace
2,
\left\lceil \frac{k-s}{rp}\right\rceil \right\rbrace & \text{if $d=2p$}\\[15pt]
\displaystyle\frac{1}{2}\binom{n}{k}\sum_{s=0}^{k-1}\binom{k}{s}\binom{n-k}{k-s}\min\left\lbrace 1+2\left\lceil\frac{\left|s-rp\right|}{2rp+r}\right\rceil , 2\left\lceil\frac{k-s}{2rp+r}\right\rceil\right\rbrace & \text{if $d=2p+1$}
    \end{cases}
\end{align*}
Unfortunately, the ceiling and the $\max/\min$ are difficult to handle, and thus it is hard to obtain a better closed formula with these tools.


\section{Acknowledgments}

\begin{sloppypar}
This work was partially supported by the LIA INFINIS/SINFIN (CNRS-CONICET-UBA, France–Argentina), the MathAm-Sud projects MATHAMSUD 20-MATH-09 and 22-MATH-02, and by Universidad Nacional de San Luis, grants PROICO 03-0723 and PROIPRO 03-2923. It was also supported by Agencia I+D+I, grants PICT 2020-00549 and PICT 2020-04064. The first author acknowledges the financial support from Lorraine University of Nancy, France, and is supported by a Ph.D. scholarship funded by Consejo Nacional de Investigaciones Científicas y Técnicas.
\end{sloppypar}


\end{document}